\documentclass[reqno]{amsart}

    \usepackage{float}
    \usepackage{amsmath}
    \usepackage{graphicx}
    \usepackage{latexsym}
    \usepackage{amsfonts}
    \usepackage{amssymb}
    \usepackage{verbatim}
    \usepackage{color}
    \theoremstyle{plain}
    \newtheorem{theorem}{Theorem}

    \newtheorem{proposition}[theorem]{Proposition}
    \theoremstyle{definition}
    \newtheorem{assumption}{Assumption}

    \newtheorem{definition}[theorem]{Definition}

    \newtheorem{remark}[theorem]{Remark}
    \newtheorem*{remark*}{Remark}

\begin{document}
    \title[Strong law of large numbers for random walks in weakly ]
    {Strong law of large numbers for random walks in weakly dependent random scenery}
    \author[Sharipov]{Sadillo Sharipov}
    \address{V.I.Romanovskiy Institute of Mathematics, Uzbekistan Academy of Sciences, Tashkent, Uzbekistan}
    \email{sadi.sharipov@yahoo.com}


\begin{abstract}
In this brief note, we study the strong law of large numbers for random walks in random scenery. Under the assumptions that the random scenery is non-stationary and satisfies weakly dependent condition with an appropriate rate, we establish strong law of large numbers for random walks in random scenery. Our results extend the known results in the literature.
\end{abstract}


    \keywords{Random walk, random scenery, local time, weak dependence, strong law of large numbers}
    \subjclass{Primary 60F15; 60G50; Secondary 62D05}
    \maketitle

\section{Introduction}
Let $\left(\Omega,\mathfrak{F},\mathbb{P} \right)$ be a probability space. The random variables we deal with are all defined on $\left(\Omega,\mathfrak{F},\mathbb{P} \right)$.
Let $\left(\xi_{i} \right)_{i\in \mathbb{Z}}$ and $\left(X_{n} \right)_{n\in \mathbb{N}}$ be two independent sequences of independent identically distributed (i.i.d.) random variables taking values in $\mathbb{R}$ and $\mathbb{Z}$, respectively. The sequence $\left(\xi_{i} \right)_{i\in \mathbb{Z}}$ is called the random scenery. The sequence $\left(X_{n} \right)_{n\in \mathbb{N}}$ is the sequence of increments of the random walk $\left(S_{n} \right)_{n\geq 0}$  defined by $S_{0}=0$ and $S_{n}:=X_{1}+\dots +X_{n}$ for $n\in \mathbb{N}$.

We define the random walk in random scenery as the process $\left(Z_{n}\right)_{n\geq 0}$ given by
\begin{equation}\label{eq1}
Z_{0}=0, \ \ Z_{n}=\sum_{k=0}^{n}\xi_{S_{k}}, \ \ n\in \mathbb{N}.
\end{equation}
In fact, $Z_{n}$ is a cumulative sum process whose summands are drawn from the scenery, the order in which the summands are drawn is determined by the path of the random walk.

A considerable number of papers have been devoted to studying the limiting behavior of \eqref{eq1}. The process \eqref{eq1} was introduced independently by Kesten and Spitzer \cite{KS} and Borodin \cite{Borodin}. For instance, Borodin \cite{Borodin} considered the case where the random walk is i.i.d., while the random scenery is independent but not necessarily identically distributed. He established a central limit theorem for \eqref{eq1} along with a Berry-Esseen bound.
At the same time, Kesten and Spitzer \cite{KS} obtained the following result: if the two sequences $\left(\xi_{i}\right)_{i \in \mathbb{Z}}$ and $\left(X_{n}\right)_{n \in \mathbb{N}}$ belong to the domain of attraction of stable laws with parameters $1 < \alpha \leq 2$ and $0 < \beta \leq 2$, respectively, then there exists $\delta > 1/2$ such that the process $\left(n^{-\delta} Z_{\lfloor nt\rfloor}\right)_{t \geq 0}$ converges in distribution in the space of continuous functions $C[0, \infty)$ to a continuous $\delta$-self-similar process with stationary increments, where $\delta = 1 - \alpha^{-1} + \left(\alpha \beta\right)^{-1}$ (here, $\lfloor a \rfloor$ is the unique integer for which $\lfloor a \rfloor \leq a < \lfloor a \rfloor+1$). They also considered the case when $0 < \alpha < 1$, with $\beta$ arbitrary, and showed that the process $\left(n^{-1/\beta} Z_{ \lfloor nt \rfloor}\right)_{t \geq 0}$ converges in distribution to a stable process with index $\beta$.
Then an intensive research has been performed to study the asymptotic behaviour of \eqref{eq1}. For a comprehensive historical overview and key results in this field, we refer the reader to the recent survey by P\`{e}ne \cite{Pene}.
From a statistical point of view, it is natural to explore relaxing the independence assumption of the random scenery in \eqref{eq1}. In this direction, Guillotin-Plantard and Prieur \cite{GNP} treated the case where the random scenery satisfies a $\theta_{2}$-weakly dependent condition with an appropriate rate and established a functional limit theorem that generalizes Kesten and Spitzer's result for the case $\alpha=\beta=2$.

To the best of our knowledge, the study of the SLLN for \eqref{eq1} has not been as extensive as for other limit theorems. The first result in this direction was established by Guillotin-Plantard \cite{GP}. She considered the dynamical $\mathbb{Z}^{d}$-random walk with a centered and independent random scenery and proved that $n^{-\tau} Z_{n} \to 0$ almost surely (a.s.) as $n \to \infty$ for every $\tau > 3/4$. The method of her proof is based on the fact that $\xi_{S_k}$ forms an associated sequence whenever the random variables $\xi_{i}, i\in \mathbb{Z}$, are independent. Later, Wang \cite{Wang2007} derived that $n^{-1}Z_{n}\to \mathbb{E}\xi_{0}$ a.s. under the more restrictive conditions that the random scenery is i.i.d.
It is worth noting that \eqref{eq1} becomes stationary when the random scenery is i.i.d. Consequently, the Birkhoff's Ergodic Theorem implies the SLLN in this case. However, $Z_{n}$ is no longer stationary if the random scenery is non-stationary. Then, it is natural to ask whether the assumptions of stationarity and independence in the random scenery can be relaxed for the SLLN. This question motivated us to study the SLLN for \eqref{eq1} in case when the random scenery is non-stationary and weakly dependent.

The purpose of this paper is to prove the validity of the SLLN for \eqref{eq1} when the random scenery is non-stationary and satisfies the $\theta_{2}$-weakly dependent condition. Furthermore, we establish the SLLN for \eqref{eq1} when the random walk exhibits long range dependence. In the special case where the scenery is identically distributed and satisfies $\theta_{1}$-weakly dependence condition, we also prove the SLLN.


Throughout this paper, the symbol $C$ denotes a positive constant which is not
necessarily the same one in each appearance and $\mathrm{I}\left(A\right)$ denotes the indicator function of
the event $A$.

Our paper is organized as follows. In Section 2, we introduce dependence settings under which we work in the
sequel. Section 3 contains main results. The proofs of our results are given in Section 4.

 \section{Assumptions and definitions}
We recall definition of the dependence coefficients which we shall use in the sequel. As a measure of dependence, we will use the weakly dependent coefficient introduced by Dedecker \cite{Dedeckeretal}.

On the Euclidean space $\mathbb{R}^m$, we define the metric
$$
d_{1}\left(x, y\right)=\sum_{i=1}^m\left|x_{i}-y_{i}\right|.
$$
Let $\Lambda=\bigcup_{m \in \mathbb{N}} \Lambda_m$, where $\Lambda_m$ is the set of Lipschitz functions $f: \mathbb{R}^{m} \rightarrow \mathbb{R}$ with respect to the metric $d_1$. If $f \in$ $\Lambda_m$, we denote by $\operatorname{Lip}(f):=\sup_{x, y} \frac{|f(x)-f(y)|}{d_{1}(x, y)}$ the Lipschitz modulus of $f$. Define the set of functions $ \widetilde{\Lambda}=\left\{f\in \Lambda: \operatorname{Lip}\left(f\right) \leq 1 \right\}$.
\begin{definition}
Let $\xi$ be $\mathbb{R}^{m}$-valued random variable defined on a probability space $\left(\Omega,\mathfrak{F},\mathbb{P} \right)$, assumed to be square integrable. For any $\sigma$-algebra $\mathcal{M}$ of $\mathfrak{F}$, we define the $\theta_{2}$-dependence coefficient
\begin{equation}\label{eq1*}
\theta_{2}\left(\mathcal{M}, \xi\right)=\sup \left\{\|\mathbb{E}\left(f(\xi) \mid \mathcal{M}\right)-\mathbb{E}(f(\xi))\|_{2}, f \in \widetilde{\Lambda}\right\},
\end{equation}
where $\|\xi \|_{2}=\left(\mathbb{E}\left|\xi \right|^{2}\right)^{1/2}$.
\end{definition}
We now define the coefficient $\theta_{k, 2}$ for a sequence of $\sigma$-algebras and a sequence of random variables.
\begin{definition} Let $\left(\xi_i\right)_{i\in \mathbb{Z}}$ be a sequence of square integrable random variables with values in $\mathbb{R}$. Let $\left(\mathcal{M}_{i}\right)_{i \in \mathbb{Z}}$ be a sequence of $\sigma$-algebras of $\mathfrak{F}$. For any $k \in \mathbb{N} \cup\{\infty\}$ and $n \in \mathbb{N}$, we define
$$
\theta_{k, 2}(n)=\max_{1 \leq l \leq k} \frac{1}{l} \sup\left\{\theta_2\left(\mathcal{M}_{p},\left(\xi_{j 1}, \ldots, \xi_{j l}\right)\right), p+n \leq j_{1}<\ldots<j_{l}\right\}
$$
and
$$
\theta_{2}\left(n\right)=\theta_{\infty, 2}\left(n\right) = \sup_{k} \theta_{k, 2}\left(n\right).
$$
The sequence $\left(\xi_i\right)_{i \in \mathbb{Z}}$ is said to be $\theta_{2}$-weakly dependent with respect to $\left(\mathcal{M}_i\right)_{i \in \mathbb{Z}}$ if $\theta_{2}(n)\to 0$ as $n \to \infty$.
\end{definition}

By substituting the norm $\| \cdot\|_{2}$ in \eqref{eq1*} with the norm $\| \cdot\|_{1}$, we obtain the $\theta_{1}$-dependence coefficient, which was first introduced by Doukhan and Louhichi \cite{DL}.

Let $\left(\mathcal{M}_{i}\right)_{i\in \mathbb{Z}}$ a sequence of $\sigma$-algebras in $\mathfrak{F}$ defined as
$\mathcal{M}_{i}=\sigma\left(\xi_{j}, j \leq i \right)$. In the following, the dependence coefficients will be defined with respect to the sequence of $\left(\mathcal{M}_{i}\right)_{i\in \mathbb{Z}}$.

Now we need notion of long range dependence in order to formulate our next result.

Let the random walk $\left(X_{n} \right)_{n\in \mathbb{N}}$ be a stationary Gaussian sequence with zero mean and correlations $r\left(i-j\right)=\mathbb{E}X_{i}X_{j}$ satisfying
\begin{equation}\label{eq2}
\sum_{i=1}^{n}\sum_{j=1}^{n}r\left(i-j\right) \sim n^{2H}L\left(n\right), \ \ n\to \infty
\end{equation}
where $0<H<1$ and $L$ is a slowly varying function at infinity. In this case, one says that $\left(X_{n}\right)_{n \in \mathbb{N}}$ exhibits long range dependence (or strong dependence).
The monograph by Beran et al. \cite{BFGK} describes in detail this notion from various points of view.

Let $i\in \mathbb{Z}$ and $n\in \mathbb{N}$. The local time $N_{n}\left(i\right)$ of the random walk $\left(S_{k}\right)_{k\geq 0}$ at point $i$ up to time $n$ is defined by
$$ N_{n}\left(i\right)=\sum_{k=0}^{n}\mathrm{I}\left(S_{k}=i\right). $$

Let $\alpha\left(n,i\right)$ be the intersection local time at the point $i\in \mathbb{Z}$ of the random walk
$\left(S_{k} \right)_{k \geq 0}$ defined by
$$ \alpha\left(n,i\right)=\sum_{k,j=0}^{n}\mathrm{I}\left(S_{k}-S_{j}=i \right). $$
Since $\sum_{i\in \mathbb{Z}}\mathrm{I}\left(S_{k}=i \right)=1$ for each $k$, and using independence assumption between the random scenery and the random walk, we infer that the process $Z_{n}$ can be rewritten as
\begin{equation}\label{eq3}
 Z_{n}=\sum_{i\in \mathbb{Z}}N_{n}\left(i\right)\xi_{i}.
\end{equation}
We need the following proposition, which provides moment bounds for the local time and intersection local time of the random walk.
\begin{proposition}
The following statements are true:
\begin{itemize}
\item[(a)] For any $p \geq 1$, there exists $C>0$ such that for all $n \geq 1$,
\begin{equation}\label{eq4}
\mathbb{E}\left(\alpha\left(n,0\right)\right)^{p}\leq C n^{3p/2};
\end{equation}
\end{itemize}
\begin{itemize}
\item[(b)] Let $\left(X_{n} \right)_{n\in \mathbb{N}}$ be a stationary Gaussian sequence with zero mean and correlations $r\left(i-j\right)=\mathbb{E}X_{i}X_{j}$ satisfying condition \eqref{eq2}. Then there exists a positive constant $C$ such that for each $n \geq 1$,
\begin{equation}\label{eq5}
\sum_{i\in \mathbb{Z}}\mathbb{E}N_{n}^{2}\left(i\right)\leq C n^{2-H}L\left(n\right).
\end{equation}
\end{itemize}
\end{proposition}
Item (a) of Proposition 3 was established by Dombry and Guillotin-Plantard \cite{DGP}, while item (b) is due to Wang \cite{Wang2003}.

\section{Statement of the Main Results}
Now we are in position to formulate our results.
\begin{theorem}\label{thm3.1}
Let $\left(X_{n} \right)_{n\in \mathbb{N}}$ be a sequence of i.i.d. random variables with $\mathbb{E}X_{0}=0$, $\mathbb{E}X_{0}^{2}=\sigma^{2}\in \left(0,\infty\right)$. Assume that the following conditions hold:\\
$(i)$ $\left(\xi_{i}^{2} \right)_{i\in \mathbb{Z}}$ is uniformly integrable;\\
$(ii)$ $\left(\xi_{i}\right)_{i\in \mathbb{Z}}$ is $\theta_{2}$-weakly dependent with $\sum_{j=0}^{\infty}\theta_{1,2}\left(j\right)<\infty$.

Then, as $n\to \infty$,
\begin{equation}\label{eq6}
\frac{Z_{n}-\mathbb{E}Z_{n}}{n} \to 0\ \  a.s.
\end{equation}
\end{theorem}
The next result deals with the validity of SLLN in the case when the random walk $\left(X_{n}\right)_{n \in \mathbb{N}}$ is long range dependent.
\begin{theorem}\label{thm3.2}
Let $\left(X_{n} \right)_{n\in \mathbb{N}}$ be a stationary Gaussian sequence with zero mean and correlations $r\left(i-j\right)$, satisfying condition \eqref{eq2}. Assume that $\left(\xi_{i}\right)_{i\in \mathbb{Z}}$ satisfies conditions $(i)$ and $(ii)$. Then, \eqref{eq6} remains valid.
\end{theorem}
The following result asserts the validity of the SLLN when the scenery is identically distributed with finite mean and exhibits $\theta_1$-weak dependence.
\begin{theorem}\label{thm3.3}
Let $\left(X_{n} \right)_{n\in \mathbb{N}}$ be a sequence of i.i.d. random variables with $\mathbb{E}X_{0}=0$, $\mathbb{E}X_{0}^{2}\in \left(0,\infty\right)$. Assume $\left(\xi_{i}\right)_{i\in \mathbb{Z}}$ is a sequence of identically distributed $\theta_{1}$-weakly dependent random variables, with $\sum_{j=0}^{\infty}\theta_{1,1}\left(j\right)<\infty$ and $\mathbb{E}\left|\xi_{0}\right|<\infty$.

Then, for each $\tau>3/4$, as $n \to \infty$,
\begin{equation}\label{eq7}
\frac{Z_{n}}{n^{\tau}} \to \mathbb{E}\xi_{0}\ \  a.s.
\end{equation}
\end{theorem}
\begin{remark}
Theorem 4 generalizes the result of Guillotin-Plantard \cite{GP} fom the specific case of $\tau=1$
with independent random scenery to the $\theta_2$-dependent condition.
Theorem 5 generalizes the corresponding result of Wang \cite{Wang2007} for i.i.d. random scenery to the non-stationary $\theta_1$-weakly dependent condition, without requiring any additional moment assumptions.
\end{remark}

\section{Proofs of the main results}
\textbf{Proof of Theorem \ref{thm3.1}.}
Due to the assumption that the random walk $\left(X_{n}\right)_{n\in \mathbb{N}}$ is i.i.d. with finite second moment, then the Hartman-Wintner law of iterated logarithm states that
\begin{equation}\label{eq7*}
 \limsup_{n \to \infty}\frac{\max_{1 \leq k \leq n}\left|S_{k} \right|}{\left(2\sigma^{2}n\log\log n\right)^{1/2}}=1 \ \ a.s.
\end{equation}
Let us pick $\delta >0$. Then we have $N_{n}\left(i\right)=0$ for all $\left|i\right|> n^{1/2+\delta}$ and
sufficiently large $n$. Therefore, taking into account \eqref{eq3},
$$
 Z_{n} = \sum_{\left|i\right|\leq \lfloor n^{\frac{1}{2}+\delta} \rfloor}N_{n}\left(i\right)\xi_{i}.
$$
Define $\widetilde{\xi}_{i}=\xi_{i}-\mathbb{E}\xi_{i}$, $i \in \mathbb{Z}$.
Let $\widetilde{\xi}_{i}^{+}=\max\left(0, \widetilde{\xi}_{i}\right)$, $\widetilde{\xi}_{i}^{-}=\max\left(0, -\widetilde{\xi}_{i}\right)$. It is known \cite{Dedeckeretal} that the variables $\widetilde{\xi}_{i}^{+}$ and $\widetilde{\xi}_{i}^{-}$ also satisfy $\theta_{2}$-weakly dependent condition. Since $\widetilde{\xi}_{i}=\widetilde{\xi}_{i}^{+}-\widetilde{\xi}_{i}^{-}$, it suffices to prove Theorem \ref{thm3.1} separately for
$\left(\widetilde{\xi}_{i}^{+}\right)_{i\in \mathbb{Z}}$ and $\left(\widetilde{\xi}_{i}^{-}\right)_{i\in \mathbb{Z}}$, thus we may and do assume that $\widetilde{\xi}_{i} \geq 0$ a.s. for all $i\in \mathbb{Z}$.

For arbitrary $\lambda> 1$ set $k_{n}=\lfloor\lambda^{n}\rfloor$, $n\in \mathbb{N}$. Here, in fact, one should throughout think of $\lambda$ as being close to 1.
For a positive integer $m\in \mathbb{N}$, there exists $n\in \mathbb{N}$ such that $k_{n}\leq m < k_{n+1}$, and $n\to \infty$ as
$m\to \infty$.\\
Since $\widetilde{\xi}_{i} \geq 0$, then it holds that
\begin{equation}\label{eq8}
\frac{Z_{m}-\mathbb{E}Z_{m}}{m} \leq
\left| \frac{Z_{k_{n+1}}-\mathbb{E}Z_{k_{n+1}}}{k_{n+1}}\right|\frac{k_{n+1}}{k_{n}}
 +\frac{\mathbb{E}Z_{k_{n+1}}-\mathbb{E}Z_{k_{n}}}{k_{n}}.
\end{equation}
From \eqref{eq8} it follows that if we have shown that
\begin{equation}\label{eq9}
\frac{Z_{k_{n}}-\mathbb{E}Z_{k_{n}}}{k_{n}} \to 0 \ \ a.s.
\end{equation}
then we would have
$$ \limsup_{m \to \infty}\left|\frac{Z_{m}-\mathbb{E}Z_{m}}{m} \right| \leq \left(\lambda-1\right) \sup_{i\in \mathbb{Z}}\mathbb{E}\xi_{i} \ \ a.s. $$
for every $\lambda > 1$ which concludes the proof.

To this end, we first estimate the variance of $Z_{k_{n}}$. By denoting $\widetilde{N}_{k_{n}}\left(i\right)=N_{k_{n}}\left(i\right)-\mathbb{E}N_{k_{n}}\left(i\right)$, we have the bounds
$$
\begin{aligned}
\operatorname{Var}\left(Z_{k_{n}}\right)
& =\mathbb{E}\left(\sum_{\left|i\right|\leq \lfloor n^{\frac{1}{2}+\delta} \rfloor}\widetilde{N}_{k_{n}}\left(i\right)
\widetilde{\xi}_{i} \right)^{2} \\
& = \sum_{\left|i\right|\leq \lfloor n^{\frac{1}{2}+\delta} \rfloor}\operatorname{Var}\left({N}_{k_{n}}\left(i\right)\right)\operatorname{Var}\left(\xi_{i}\right)+\sum_{\left|i\right|\leq \lfloor n^{\frac{1}{2}+\varepsilon} \rfloor}
\sum_{\left|j\right|\leq \lfloor n^{\frac{1}{2}+\delta} \rfloor}\mathbb{E}\widetilde{N}_{k_{n}}\left(i\right)\widetilde{N}_{k_{n}}\left(j\right)\mathbb{E} \widetilde{\xi}_{i} \widetilde{\xi}_{j} \\
& \leq \sum_{\left|i\right|\leq \lfloor n^{\frac{1}{2}+\delta} \rfloor}\operatorname{Var}\left(N_{k_{n}}\left(i\right)\right)\operatorname{Var}\left(\xi_{i}\right)+\sum_{\left|i\right|\leq \lfloor n^{\frac{1}{2}+\varepsilon} \rfloor}
\sum_{\left|j\right|\leq \lfloor n^{\frac{1}{2}+\delta} \rfloor}\left|\mathbb{E}\widetilde{N}_{k_{n}}\left(i\right)\widetilde{N}_{n_{k}}\left(j\right)\right|\left|\mathbb{E} \widetilde{\xi}_{i} \widetilde{\xi}_{j}\right| \\
& \leq \sum_{\left|i\right|\leq \lfloor n^{\frac{1}{2}+\delta} \rfloor}\operatorname{ Var}\left(N_{k_{n}}\left(i\right)\right)\operatorname{Var}\left(\xi_{i}\right)+\sum_{\left|i\right|\leq \lfloor n^{\frac{1}{2}+\delta} \rfloor}
\sum_{\left|j\right|\leq \lfloor n^{\frac{1}{2}+\delta} \rfloor}\|\widetilde{N}_{k_{n}}\left(i\right)\|_{2}\|\widetilde{N}_{k_{n}}\left(j\right)\|_{2}\left|\mathbb{E} \widetilde{\xi}_{i} \widetilde{\xi}_{j}\right| \\
& \leq \sum_{\left|i\right|\leq \lfloor n^{\frac{1}{2}+\delta} \rfloor}\operatorname{Var}\left(N_{k_{n}}\left(i\right)\right)\operatorname{Var}\left(\xi_{i}\right)+
\sum_{\left|i\right|\leq \lfloor n^{\frac{1}{2}+\delta} \rfloor}\operatorname{Var}\left(N_{k_{n}}\left(i\right)\right)
\sum_{\left|j\right|\leq \lfloor n^{\frac{1}{2}+\delta} \rfloor, j \neq i}\left|\mathbb{E} \widetilde{\xi}_{i} \widetilde{\xi}_{j}\right|,
\end{aligned}
$$
where in the last step we used that $2ab\leq a^{2}+b^{2}$ for any real numbers $a,b$.

We now have to control the covariance $\left|\mathbb{E}\widetilde{\xi}_{i}\widetilde{\xi}_{j}\right|$. Using the Cauchy-Schwarz inequality, it is easy to see that for $i<j$,
\begin{equation}\label{eq10}
\begin{aligned}
\left|\mathbb{E}\widetilde{\xi}_{i}\widetilde{\xi}_{j}\right|
& = \left|\mathbb{E}\left(\widetilde{\xi}_{i}\mathbb{E}\left(\widetilde{\xi}_{j}| \mathcal{M}_{i}\right)\right)\right| \\
& \leq\left(\mathbb{E}\left|\widetilde{\xi}_{i}\right|^{2} \right)^{1/2}\left(\mathbb{E}\left(\mathbb{E}\left(\widetilde{\xi}_{j}|\mathcal{M}_{i}\right)\right)^{2}\right)^{1/2} \\
& \leq\left(\operatorname{Var}\left(\widetilde{\xi}_{i}\right) \right)^{1/2} \theta_{1,2}\left(j-i\right).
\end{aligned}
\end{equation}
Substituting the above bounds and from conditions $(i)$ and $(ii)$, we deduce
\begin{equation}\label{eq11}
\operatorname{Var}\left(Z_{k_{n}}\right) \leq C \left(\sup_{i\in \mathbb{Z}}\operatorname{Var}\left(\xi_{i}\right)+\sqrt{\sup_{i\in \mathbb{Z}}\operatorname{Var}\left(\xi_{i}\right)}\sum_{l=0}^{\infty}\theta_{1,2}\left(l\right) \right)\sum_{\left|i\right|\leq \lfloor n^{\frac{1}{2}+\delta} \rfloor}\operatorname{Var}\left(N_{k_{n}}\left(i\right)\right).
\end{equation}

Consequently, by the Chebyshev inequality and \eqref{eq4}, for any $\varepsilon>0$, we infer that
$$
\begin{aligned}
\sum_{n=1}^{\infty}\mathbb{P}\left(\frac{Z_{k_{n}}-\mathbb{E}Z_{k_{n}}}{k_{n}}>\varepsilon\right)
& \leq \sum_{n=1}^{\infty}\frac{1}{k_{n}^{2}}\operatorname{Var}\left(Z_{k_{n}}\right) \\
& \leq C \sum_{n=1}^{\infty}\frac{1}{k_{n}^{2}} \sum_{\left|i\right|\leq \lfloor n^{\frac{1}{2}+\delta} \rfloor}\operatorname{Var}\left(N_{k_{n}}\left(i\right)\right) \\
&\leq C \sum_{n=1}^{\infty}\frac{k_{n}^{3/2}}{k_{n}^{2}} <\infty,
\end{aligned}
$$
which proves \eqref{eq9} and hence \eqref{eq6}. Theorem is proved.

\textbf{Proof of Theorem \ref{thm3.2}.}
We proceed as in the proof of Theorem \ref{thm3.1}. First, taking into account that the random walk $\left(X_{n} \right)_{n\in \mathbb{N}}$ is long range dependent, then by law of iterated logarithm, we deduce
$  S_{n}=O\left(n^{H}\sqrt{L\left(n\right)}\log n\right)$ a.s.
Hence, for all $\left|i \right|> n^{H+\delta}$, where $\delta$ is an arbitrary positive number, it yields $N_{n}\left(i\right)=0$.
This implies that
$$  Z_{n} = \sum_{\left|i\right|\leq \lfloor n^{H+\delta} \rfloor}N_{n}\left(i\right)\xi_{i}. $$
Let $\lambda>1$ and set $k_{n}=\lfloor \lambda^{n} \rfloor$.  According to \eqref{eq8} we have to prove \eqref{eq9}. Note that \eqref{eq10} and \eqref{eq11} also remains valid.
Thus, it remains to prove the summability of $\mathbb{P}\left(\frac{Z_{k_{n}}-\mathbb{E}Z_{k_{n}}}{k_{n}} \right)$.
Using the Chebyshev inequality and then \eqref{eq5}, we have for each $\varepsilon > 0$,
$$
\begin{aligned}
\sum_{n=1}^{\infty}\mathbb{P}\left(\frac{Z_{k_{n}}-\mathbb{E}Z_{k_{n}}}{k_{n}}>\varepsilon\right)
& \leq \sum_{n=1}^{\infty}\frac{1}{k_{n}^{2}}\operatorname{Var}\left(Z_{k_{n}}\right) \\
& \leq C \sum_{n=1}^{\infty}\frac{1}{k_{n}^{2}} \sum_{\left|i\right|\leq \lfloor n^{H+\delta} \rfloor}\operatorname{Var}\left(N_{k_{n}}\left(i\right)\right) \\
& \leq C \sum_{n=1}^{\infty}\frac{k_{n}^{2-H}L\left(k_{n}\right)}{k_{n}^{2}} <\infty.
\end{aligned}
$$
Hence, by the Borel--Cantelli lemma we arrive at \eqref{eq9}, which due to \eqref{eq8} concludes the proof of Theorem \ref{thm3.2}.

\textbf{Proof of Theorem \ref{thm3.3}.} Without restricting the generality, we can assume $\xi_{i}\geq 0$ a.s. for all $i \in \mathbb{Z}$.
For each $n \in \mathbb{N}$ and $i\in \mathbb{Z}$, we define the truncated random variables $\zeta_{i}=\xi_{i}\mathrm{I}\left(\xi_{i}<n\right)$. It is known \cite{Dedeckeretal} that the random variables $\zeta_{i}$, $i\in \mathbb{Z}$ satisfy $\theta$-weakly dependent condition and for any fixed path of the random walk, $\theta_{k,1}$-weakly dependent coefficients of the sequences $\zeta=\left(\zeta_{i}\right)_{i\in \mathbb{Z}}$ and $\xi=\left(\xi_{i}\right)_{i\in \mathbb{Z}}$ satisfies $\theta_{k,1}^{\zeta}\left(n\right) \leq \theta_{k,1}^{\xi}\left(n\right)$, $n,k \geq 1$.
For any $\lambda > 1$, set $k_{n}=\lfloor \lambda^{n}\rfloor$.

First, we show that as $n \to \infty$
\begin{equation}\label{eq12}
\frac{Z_{k_{n}}}{k_{n}^{\tau}} \to 0 \ \ a.s.
\end{equation}
To this end, we note that \eqref{eq7*} holds and thus we denote
$$
\widetilde{Z}_{k_{n}}=\sum_{\left|i\right|\leq \lfloor n^{\frac{1}{2}+\delta} \rfloor}N_{k_{n}}\left(i\right)\zeta_{i}.
$$
Observe that for all $i \in \mathbb{Z}$, it holds
\begin{equation}\label{eq12*}
\sum_{n=1}^{\infty}\mathbb{P}\left(N_{k_{n}}\left(i\right)\xi_{i} \neq N_{k_{n}}\left(i\right)\zeta_{i} \right) = \sum_{n=1}^{\infty}\mathbb{P}\left(\xi_{i} \geq n \right)=\sum_{n=1}^{\infty}\mathbb{P}\left(\xi_{0} \geq n \right) < \infty.
\end{equation}
Hence, $\mathbb{P}\left(N_{k_{n}}\left(i\right)\xi_{i} \neq N_{k_{n}}\left(i\right)\zeta_{i}, i.o. \right)=0$ by the Borel--Cantelli lemma.
Thus, \eqref{eq12} holds if and only if
\begin{equation}\label{eq13}
\frac{\widetilde{Z}_{k_{n}}}{k_{n}^{\tau}} \to 0 \ \ a.s.
\end{equation}
From the previous arguments made in the proof of Theorem 4, it follows that
$$ \left|\operatorname{cov} \left(\zeta_{i}, \zeta_{j}\right)\right|  \leq \mathbb{E}\left|\zeta_{0}-\mathbb{E}\zeta_{0} \right| \theta_{1,1}^{\zeta}\left(j-i\right). $$
Then
$$
\begin{aligned}
\mathbb{E}\left(\widetilde{Z}_{k_{n}}\right)^{2}
& = \mathbb{E}\left(\sum_{\left|i\right|\leq \lfloor n^{\frac{1}{2}+\delta} \rfloor}N_{k_{n}}
\left(i\right)\zeta_{i} \right)^{2} \\
& \leq \sum_{\left|i\right|\leq \lfloor n^{\frac{1}{2}+\delta} \rfloor}\mathbb{E}N_{k_{n}}^{2}\left(i\right)\mathbb{E}\zeta_{i}^{2}+\sum_{\left|i\right|\leq \lfloor n^{\frac{1}{2}+\delta} \rfloor}
\sum_{\left|j\right|\leq \lfloor n^{\frac{1}{2}+\delta} \rfloor}\left|\mathbb{E}N_{k_{n}}\left(i\right)N_{k_{n}}\left(j\right)\right|\left|\operatorname{cov}\left(\zeta_{i}, \zeta_{j}\right)\right| \\
& \leq \sum_{\left|i\right|\leq \lfloor n^{\frac{1}{2}+\delta} \rfloor}\mathbb{E}N_{k_{n}}^{2}\left(i\right)\mathbb{E}\zeta_{i}^{2}+
\sum_{\left|i\right|\leq \lfloor n^{\frac{1}{2}+\delta} \rfloor}\mathbb{E}N_{k_{n}}^{2}\left(i\right)
\sum_{\left|j\right|\leq \lfloor n^{\frac{1}{2}+\delta} \rfloor, j \neq i}\left|\operatorname{cov}\left(\zeta_{i}, \zeta_{j}\right)\right| \\
& \leq C \sum_{\left|i\right|\leq \lfloor n^{\frac{1}{2}+\delta} \rfloor}\mathbb{E}N_{n_{k}}^{2}\left(i\right),
\end{aligned}
$$
where $C=\mathbb{E}\zeta_{0}^{2}+\left(\mathbb{E}\zeta_{0}^{2}\right)^{1/2}\sum_{l=0}^ {\infty}\theta_{1,1}^{\zeta}\left(l\right)$.

Further, from \eqref{eq4}, it follows that for all $\varepsilon>0$,
$$
\begin{aligned}
\sum_{n=1}^{\infty}\mathbb{P}\left(\frac{1}{k_{n}^{\tau}}\widetilde{Z}_{k_{n}} >\varepsilon\right)
& \leq C \sum_{n=1}^{\infty}\frac{1}{k_{n}^{2\tau}}\mathbb{E} \left(\widetilde{Z}_{k_{n}}\right)^{2} \\
& \leq C \sum_{n=1}^{\infty}\frac{1}{k_{n}^{2\tau}} \mathbb{E}\xi_{0}^{2}\mathrm{I}\left(\xi_{0}<n\right)\sum_{\left|i\right|\leq \lfloor n^{\frac{1}{2}+\varepsilon} \rfloor}\mathbb{E}N_{k_{n}}^{2}\left(i\right) \\
& \leq C \sum_{n=1}^{\infty}\frac{1}{n^{2}} \mathbb{E}\xi_{0}^{2}\mathrm{I}\left(\xi_{0}<n\right),
\end{aligned}
$$
Observe that
$$
\begin{aligned}
\sum_{n=1}^{\infty}\frac{1}{n^{2}} \mathbb{E}\xi_{0}^{2}\mathrm{I}\left(\xi_{0}<n\right)
& =\sum_{n=1}^{\infty}\frac{1}{n^{2}}\sum_{j=1}^{n} \mathbb{E}\xi_{0}^{2}\mathrm{I}\left(j-1 \leq \xi_{0}< j \right) \\
& \leq \sum_{j=1}^{\infty} j^{2} \mathbb{P}\left(j-1 \leq \xi_{0}< j \right) \sum_{n=j}^{\infty}\frac{1}{n^{2}} \\
& \leq C\sum_{j=1}^{\infty} j \mathbb{P}\left(j-1 \leq \xi_{0}< j \right) \\
& = C \sum_{j=1}^{\infty}\mathbb{P}\left(\xi_{0} \geq j \right)< \infty.
\end{aligned}
$$
From the Borel--Cantelli lemma, we obtain \eqref{eq13}. Combining \eqref{eq12*} with \eqref{eq13}, we arrive at \eqref{eq12}.

To finish the proof, it remains to note that for any $n \in \mathbb{N}$ there exists $k_{n} \in \mathbb{N}$ such that $k_{n-1}\leq n < k_{n}$. Hence, by the monotonicity of $Z_{n}$,
$$ \frac{1}{\lambda^{\tau}} \frac{Z_{k_{n-1}}}{k_{n-1}^{\tau }} \leq \frac{Z_{n}}{n^{\tau }} \leq \lambda^{\tau} \frac{Z_{k_{n}}}{k_{n}^{\tau}} $$
which follows that
$$ \frac{1}{\lambda^{\tau}}\mathbb{E}\xi_{0} \leq \liminf_{n\to \infty} \frac{Z_{n}}{n^{\tau}} \leq \limsup_{n\to \infty}\frac{Z_{n}}{n^{\tau}} \leq \lambda^{\tau} \mathbb{E}\xi_{0}. $$
Since $\lambda$ may be arbitrarily close to 1, we get \eqref{eq7}. Theorem \ref{thm3.3} is proved.


\begin{thebibliography}{99}

\bibitem{BFGK} Beran, J., Feng, Y., Ghosh, S. and Kulik, R.
\newblock
\newblock\emph{Long-Memory Processes: Probabilistic Properties and Statistical Methods}. Springer, Berlin Heidelberg, 2013.


\bibitem{Borodin} Borodin, A.N.
\newblock Limit theorems for sums of independent random variables defined on a transient random walk. Investigations in the theory of
probability distributions, IV.
\newblock \emph{Zap. Nauchn. Sem. Leningrad. Otdel. Mat. Inst. Steklov} \textbf{85}: 17--29, 237, 244, 1979.

\bibitem{Dedeckeretal} Dedecker, J., Doukhan, P., Lang, G., Le\'{o}n, J. R.,  Louhichi, S. and Prieur, C.
\newblock
\newblock\emph{Weak Dependence: With Examples and Applications. Lect. Notes in Stat.}, \textbf{190}. Springer, New-York, 2007.

\bibitem{DGP} Dombry, C. and Guillotin-Plantard, N.
\newblock Discrete approximation of a stable self-similar stationary increments process.
\newblock \emph{Bernoulli} \textbf{15}: 195--222, 2009.

\bibitem{DL} Doukhan, P. and Louhichi, S.
\newblock A new weak dependence condition and applications to moment inequalities.
\newblock \emph{Stochastic Process. Appl} \textbf{84}: 313--342, 1999.


\bibitem{GP} Guillotin-Plantard, N.
\newblock Dynamic $\mathbb{Z}^{d}$-Random Walks in a Random Scenery:
A Strong Law of Large Numbers.
\newblock \emph{Journal of Theoretical Probability} \textbf{14}:241--260, 2001.

\bibitem{GNP} Guillotin-Plantard, N. and Prieur, C.
\newblock Limit theorem for random walk in weakly dependent
random scenery.
\newblock \emph{Annales de l'Institut Henri Poincare - Probabilites et Statistiques} \textbf{46}:1178--1194, 2010.

\bibitem{KS} Kesten, H. and Spitzer, F.
\newblock A limit theorem related to a new class of self-similar processes.
\newblock \emph{Z. Wahrsch. Verw. Gebiete} \textbf{50}: 5--25, 1979.

\bibitem{Pene} P\`{e}ne, F.
\newblock Random walks in random sceneries and related models.
\newblock \emph{ESAIM: Proceedings and Surveys} {\bf 68}: 35--51, 2020.

\bibitem{Wang2003} Wang, W.
\newblock Weak convergence to fractional Brownian motion in Brownian motion
\newblock \emph{Journal of Theoretical Probability} {\bf 126}: 203--220, 2003.

\bibitem{Wang2007} Wang, W.
\newblock Strong laws of large numbers for random walks in
random sceneries
\newblock \emph{Acta Mathematicae Applicatae Sinica, English Series} {\bf 23}: 495--500, 2007.

\end{thebibliography}
\end{document}